\documentclass[preprintnumbers,amsmath,amssymbm,preprint]{revtex4}
\usepackage{epsfig}
\usepackage{graphicx}

\begin{document}
\title{Sisyphus random walks in the presence of moving traps}
\author{Shahar Hod}
\affiliation{The Ruppin Academic Center, Emeq Hefer 40250, Israel}
\affiliation{ }
\affiliation{The Hadassah Institute, Jerusalem 91010, Israel}
\date{\today}

\begin{abstract}
\ \ \ It has recently been proved that, in the presence of a static
absorbing trap, Sisyphus random walkers with a restart mechanism are
characterized by {\it exponentially} decreasing asymptotic survival
probability functions. Interestingly, in the present compact paper
we prove analytically that, in the presence of a moving trap whose
velocity approaches zero asymptotically in time as $v_{\text{trap}}\sim
1/t$, the survival probabilities of the Sisyphus walkers are
dramatically changed into inverse {\it power-law} decaying tails. 
\end{abstract}
\bigskip
\maketitle

\section{Introduction}

The dynamics of random walkers in the presence of absorbing traps
have attracted much attention from physicists and mathematicians
during the last three decades \cite{RW1,RW2,RW3,RW4,RW5,RW6}. In
particular, random walk models provide useful physical insights in
various studies of complex adaptive systems \cite{DV2,Hodetl}, granular
segregation models \cite{DV3}, and polymer adsorption phenomena
\cite{DV4}. Random walk models are also useful in the mathematical
analysis of absorbing-state phase transitions \cite{DV1}. In
addition, some biological toy models of epidemic spreading are based
on the physical and mathematical properties of biased random walkers \cite{DV5}.

One of the most important mathematical quantities that 
characterizes the dynamics of random walkers in the presence of
absorbing traps is the survival probability function $S(t)$ of the
population. This physically interesting quantity describes the
time-dependent probability that a random walker has not been
absorbed by the trap after taking $t$ steps. It is well known that,
for biased random walkers who are characterized by a constant drift
velocity towards the absorbing trap, the asymptotic survival
probability decays {\it exponentially} in time \cite{DV2,Jsp},
\begin{equation}\label{Eq1}
S(t\to\infty)\sim t^{-3/2}e^{-\gamma t}\ \ \ \ \text{with}\ \ \ \
\gamma(q)=2\sqrt{q(1-q)}\  ,
\end{equation}
where $q>1/2$ is the probability of the random walker to step
towards the trap (and $1-q<1/2$ is the complementary probability of
the walker to step in the opposite direction of the trap).

Interestingly, it has been proved in \cite{Jsp} that if the
absorbing trap moves with a {\it constant} (time-independent)
velocity which equals the average drift velocity $v_{\text{drift}}=2q-1>0$
of the walkers, then the temporal behavior of the survival
probability function changes dramatically. In particular, one finds
that in this case the survival probability decays asymptotically as
an inverse {\it power} of time (instead of exponentially in time)
\cite{DV2,Jsp}:
\begin{equation}\label{Eq2}
S(t)\sim t^{-1/2}\ \ \ \ \text{for}\ \ \ \ v_{\text{trap}}=2q-1\  .
\end{equation}

In recent years there is a growing interest of physicists and mathematicians in analyzing the properties of the Sisyphus random
walk model \cite{RS1,RS2,RS3,RS4,RS5,Sisy,Hodanp,Hodbia},
a stochastic process with a simple restart mechanism in which the
random walker, regardless of her current position, has a finite
probability $1-q$ to return to her initial position (which, we
assume, is located at the origin). In this physically interesting model, the
asymptotic survival probability of the walkers in the presence of a
static trap decays {\it exponentially} in time for all physically
acceptable values $q\in(0,1)$ of the biased jumping probability
\cite{Hodbia}:
\begin{equation}\label{Eq3}
S(t\to\infty)\sim e^{-\gamma_{\text{Sis}} t}\ \ \ \ \text{with}\ \ \
\ \gamma_{\text{Sis}}(x_T\gg1)=-\ln[1-q^{x_T}(1-q)]\  ,
\end{equation}
where $x_T$ is the location of the trap (the initial distance
between the walkers and the absorbing trap).

Similarly to the case of ordinary random walks (with no restart mechanisms) \cite{Jsp}, it is
natural to expect that the asymptotic late-time behavior of the
survival probability function of the Sisyphus random walkers may be
changed dramatically if the absorbing trap moves away from the
starting point of the walkers. Based on this expectation, we here
raise the following physically interesting question: How fast should
the trap move in order to change the qualitative behavior of the
survival probability function of the Sisyphus random walkers from an
{\it exponentially} decaying function [see Eq. (\ref{Eq3})] into an inverse
{\it power-law} decaying function?

As we shall explicitly prove below, the answer to the above stated
question is highly non-trivial. In particular, in this paper we
shall reveal the somewhat surprising fact that, contrary to the
familiar case of ordinary random walk models with no restart
mechanisms in which the survival probability function decays as an
inverse power of time if the trap moves with a {\it constant}
(non-decreasing in time) velocity [see Eq. (\ref{Eq2})] \cite{Jsp},
the survival probability function that characterizes the dynamics of
the Sisyphus random walkers decays asymptotically as an inverse
power of time in the presence of moving traps whose velocities go to
{\it zero} asymptotically as $v_{\text{trap}}\sim 1/t$.

\section{Description of the system}

We analyze the survival probability function $S[t;x_T(t),q]$
which characterizes the dynamics of biased Sisyphus random walkers
in the presence of a moving trap. A Sisyphus agent, who starts her
dynamics at the origin $x_0\equiv x(t=0)=0$, is characterized by a finite 
probability $q$ to take a unit step towards the
absorbing trap, whose location $x_T(t)>0$ is a non-decreasing
function of time, and a complementary probability
$1-q$ to jump back all the way to her initial location at the
origin.

The dynamics of a biased Sisyphus random walker is therefore
governed by the mathematically compact jumping rule
\begin{equation}\label{Eq4}
x(t+1)=
\begin{cases}
x(t)+1& \text{with probability}\ \ q\ ;
\\ 0 & \text{with probability}\ \ 1-q\  ,
\end{cases}
\end{equation}
supplemented by the normalized boundary condition
\begin{equation}\label{Eq5}
N_{\text{tot}}(t)=1 \ \ \ \ \text{for}\ \ \ \ t<x_T(t)\ .
\end{equation}
Here
\begin{equation}\label{Eq6}
N_{\text{tot}}(t)=\sum_{k=0}^{k=x_T(t)-1}N_k(t)\
\end{equation}
is the time-dependent total number of Sisyphus walkers [normalized
to $N_{\text{tot}}(t=0)=1$ at the beginning of the dynamics] who
have not been absorbed by the moving trap after taking $t$ steps
according to the jumping rule (\ref{Eq4}), and $N_k(t)$ is the number
of agents who are located at $0\leq x(t)=k<x_T(t)$ at time $t$.

\section{The characteristic recurrence relation for the survival probability
function}

In the present section we shall use analytical techniques in order
to derive the recurrence relation [see Eq. (\ref{Eq12}) below] that 
determines the functional behavior of the survival probability
function $S(t)$, the time-dependent (monotonically decreasing)
fraction of Sisyphus walkers who have not been absorbed by the
moving trap after taking $t$ steps. This important physical quantity
is defined by the relation [see the normalization (\ref{Eq5})]
\begin{equation}\label{Eq7}
S[t;x_T(t),q]\equiv N_{\text{tot}}(t)\  .
\end{equation}

We first note that, taking cognizance of the jumping rule
(\ref{Eq4}), one finds that the number
\begin{equation}\label{Eq8}
d(t)\equiv N_{\text{tot}}(t-1)-N_{\text{tot}}(t)\
\end{equation}
of walkers who, after taking exactly $t$ steps, are absorbed by the
trap at $x=x_T(t)$ is given by the relation
\begin{equation}\label{Eq9}
d(t)=q\cdot N_{x_T(t)-1}(t-1)=q^2\cdot
N_{x_T(t)-2}(t-2)=\cdots =Q\cdot N_0[t-x_T(t)]\  ,
\end{equation}
where
\begin{equation}\label{Eq10}
Q\equiv q^{x_T(t)}\  .
\end{equation}
Here we have used the jumping rule (\ref{Eq4}), according to which
an agent who is located at the origin at time $t-x_T(t)$ has to go
$x_T(t)$ steps in a row to the right (towards the absorbing trap) in
order to be absorbed at time $t$ by the trap, which is now located
at $x_T(t)$ [remember that each step to the right is characterized
by the finite probability $q$].

The normalized number $N_0[t-x_T(t)]$ of walkers who are located at
the origin at time $t-x_T(t)$ is given by
\begin{equation}\label{Eq11}
N_0[t-x_T(t)]=(1-q)\cdot N_{\text{tot}}[t-x_T(t)-1]\
.
\end{equation}
Here we have used the jumping rule (\ref{Eq4}), according to which a
Sisyphus random walker has a finite probability $1-q$ to
jump all the way back to her initial position at the origin. Taking
cognizance of Eqs. (\ref{Eq7}), (\ref{Eq8}), (\ref{Eq9}), (\ref{Eq10}), and
(\ref{Eq11}), one obtains the recurrence relation
\begin{equation}\label{Eq12}
S(t)=S(t-1)-Q\cdot(1-q)\cdot S[t-x_T(t)-1]\  ,
\end{equation}
which determines the time-dependent fraction of the Sisyphus walkers
who have not been absorbed by the moving trap until the $t$-th time
step.

It is worth emphasizing the fact that the functional
behavior of the survival probability function $S(t)$, which is governed by
the recurrence relation (\ref{Eq12}), depends on the temporal
behavior of the function $x_T(t)$, the location of the moving trap.

In particular, as emphasized above, it has recently been proved
\cite{Hodbia} that, for a fixed location $x_T$ of the absorbing trap, 
the Sisyphus random walkers are characterized
by {\it exponentially} decaying asymptotic survival probability
functions of the form (\ref{Eq3}). 

Intriguingly, as we shall prove
explicitly in the next section, the qualitative behavior of the
late-time survival probability function of the Sisyphus random
walkers may be changed dramatically into an inverse {\it power-law}
tail of the form
\begin{equation}\label{Eq13}
S(t\to\infty)={{\alpha}\over{t^{\beta}}}\
\end{equation}
if the absorbing trap moves to the right (away from the initial
position of the walkers) with an asymptotically {\it vanishing}
velocity of the form $v_{\text{trap}}(t\to\infty)\sim 1/t$.

\section{Survival probabilities in the Sisyphus random walk model with moving traps}

In the present section we shall analyze the survival probabilities of Sisyphus random walkers 
in the presence of moving traps [$x_T=x_T(t)$]. 
Substituting the ansatz (\ref{Eq13}) for the
asymptotic functional behavior of the time-dependent survival
probability function into the analytically derived recurrence relation (\ref{Eq12}), one
finds the characteristic equation
\begin{equation}\label{Eq14}
1=\Big(1-{{1}\over{t}}\Big)^{-\beta}-q^{x_T(t)}\cdot(1-q)\cdot
\Big[1-{{x_T(t)+1}\over{t}}\Big]^{-\beta}\
\end{equation}
for the time-dependent location $x_T(t)$ of the moving absorbing
trap.

Using the asymptotic relation
\begin{equation}\label{Eq15}
\Big(1-{{c}\over{t}}\Big)^d=1-{{c\cdot d}\over{t}}+O(t^{-2})\ \ \ \
\text{for}\ \ \ \ t\to\infty\  ,
\end{equation}
one can express (\ref{Eq14}) in the form \cite{Noteund}
\begin{equation}\label{Eq16}
q^{x_T(t)}={{\beta}\over{1-q}}\cdot {{1}\over{t}}+O(t^{-2})\  ,
\end{equation}
which yields the simple temporal behavior \cite{Noteflo}
\begin{equation}\label{Eq17}
x_T(t)=a\cdot\ln t+b\ \ \ \ \text{for}\ \ \ \ t\to\infty\  
\end{equation}
for the location of the moving trap, where
\begin{equation}\label{Eq18}
a={{1}\over{\ln(1/q)}}\ \ \ \ ; \ \ \ \
b={{\ln\big({{\beta}\over{1-q}}\big)}\over{\ln(q)}}\  .
\end{equation}
Interestingly, from Eq. (\ref{Eq17}) one finds that the velocity of
the moving trap {\it vanishes} asymptotically:
\begin{equation}\label{Eq19}
v_{\text{trap}}\equiv {{dx_T(t)}\over{dt}}={{a}\over{t}}\to0\ \ \ \
\text{for}\ \ \ \ t\to\infty\  .
\end{equation}

\section{Summary and Discussion}

In the present paper we have used analytical techniques in order to
analyze the late-time functional behavior of the survival
probability function $S(t)$ which characterizes the dynamics of
biased Sisyphus random walkers with the restart mechanism (\ref{Eq4}) in the presence of absorbing traps.

It is well known that the survival probability function of ordinary
random walkers (who have no restart mechanisms) decays exponentially
in time [see Eq. (\ref{Eq1})] \cite{DV2,Jsp} if the jumping
probability $q$ towards the trap is larger than $1/2$. On the other
hand, if the absorbing trap moves away from the walkers with a {\it
constant} (non-decreasing in time) velocity which equals the drift
velocity of the biased walkers,
$v_{\text{trap}}=v_{\text{drift}}=2q-1>0$, then the qualitative
behavior of the survival function changes dramatically \cite{Jsp} -- it now
decays asymptotically slower, as an inverse {\it power} of time [see
Eq. (\ref{Eq2})] \cite{Jsp}.

As for the case of Sisyphus random walkers with fixed absorbing
traps, it has recently been proved in \cite{Hodbia} that the
asymptotic survival probability function of the walkers decays {\it
exponentially} in time for all values $0<q<1$ of the biased jumping
probability [see Eq. (\ref{Eq3})]. Motivated by the interesting
results of \cite{Jsp} for the case of moving traps in ordinary
random walk models (with no restart mechanisms), in the present
paper we have raised the following physically intriguing question:
How fast should the trap move in the Sisyphus random walk model in
order to change the qualitative temporal behavior of the late-time
survival probability from an exponentially decaying function into an
inverse power-law decaying function?

Using analytical techniques, we have revealed the fact that the
answer to the above stated question is highly non-trivial. In
particular, contrary to the constant velocity of the trap which is
needed in the ordinary random walk model in order to change the
behavior of its survival probability function into an inverse
power-law decaying tail , in the Sisyphus random walk model the
required velocity of the trap is remarkably {\it smaller}. In
particular, we have explicitly proved that if the trap moves
according to the simple law $x_T(t)=a\cdot\ln t+b$ (that is, with an
asymptotically {\it vanishing} velocity, $v_{\text{trap}}=a/t\to
0$), then the survival probability function of the Sisyphus walkers
decays asymptotically as an inverse power of time [see Eqs.
(\ref{Eq13}), (\ref{Eq17}), and (\ref{Eq18})] \cite{Notealph,Noteosc}:
\begin{equation}\label{Eq20}
S(t\to\infty)={{\alpha}\over{t^{\beta}}}\ \ \ \ \text{with}\ \ \ \
\beta=(1-q)\cdot q^b\  .
\end{equation}

Finally, it is interesting to point out that the power-law index
$\beta$ of the survival probability function (\ref{Eq20}) in the Sisyphus random walk model is universal
in the sense that it is independent of the explicit value of the
parameter $a$ in the time-dependent velocity profile (\ref{Eq19}) of
the absorbing trap \cite{Notean0}.

\newpage
\noindent
{\bf ACKNOWLEDGMENTS}
\bigskip

This research is supported by the Carmel Science Foundation. I would
like to thank Yael Oren, Arbel M. Ongo, Ayelet B. Lata, and Alona B.
Tea for helpful discussions.



\begin{thebibliography}{99}
\bibitem{Noteund} We assume that the velocity of the moving trap
vanishes asymptotically [see Eq. (\ref{Eq19}) below], in which case
$x_T(t)/t\to0$ in the late-time asymptotic limit $t\to\infty$.
\bibitem{Noteflo} We shall henceforth assume that $x_T(t)$ stands
for $\left\lfloor{x_T(t)}\right\rfloor$, where
$\left\lfloor{x_T}\right\rfloor$ is the floor function (the largest
integer which is smaller than or equal to $x_T$). Since our main
focus in this paper is on the asymptotic large-$t$ behavior of the
survival probability function, then
$x_T(t)/\left\lfloor{x_T(t)}\right\rfloor=1+O(1/\ln t)$
asymptotically [see Eq. (\ref{Eq17})].
\bibitem{Notealph} The pre-factor $\alpha$ in Eq. (\ref{Eq13}) can be estimated from
the simple probability relation $S[t=x_T(t)]=1-q^{x_T(t)}$.
\bibitem{Noteosc} In order to verify the validity of our analytical treatment, 
one may substitute the exact functional relation $x_T(t)=a\cdot\ln t+b$ 
for the time-dependent location of the absorbing trap into Eq. (\ref{Eq12}) which 
determines the survival probability function of the system. 
In this case one finds the time-dependent functional relation 
$S(t)={{\alpha}\over{t^{\beta}}}\cdot[1+\epsilon\cdot{{\ln t}\over{t}}+O(t^{-1})]$ 
for the survival probability of the Sisyphus random walkers in the presence of the moving trap, 
where $\beta=q^b(1-q)$ [see Eq. (\ref{Eq18})] and $\epsilon={{[q^b(1-q)]^2}\over{\ln(1/q)}}$. 
This result implies that the asymptotic relation (\ref{Eq20}) is valid in 
the $t\gg\ln t$ regime.
\end{thebibliography}
\end{document}